\documentclass{theclass}
\usepackage[T1]{fontenc}
\usepackage[utf8]{inputenc}
\usepackage[british]{babel}
\usepackage{mathtools}
\usepackage{microtype}
\usepackage{tikz}
\newcommand{\gi}{\gamma_{\mathrm{i}}}
\newcommand{\ggr}{\gamma_{\mathrm{gr}}}
\newcommand{\HH}{H}
\newcommand{\TT}{\mathcal T}

\newtheoremstyle{openquestion}{\topsep}{\topsep}{}{}{\bfseries}{:}{.5em}
  {\thmname{#1}\thmnumber{ #2}\thmnote{ \normalfont #3}}
\theoremstyle{openquestion}
\newtheorem{question}[thm]{Question}

\hypersetup{colorlinks=true,linkcolor=red,citecolor=green,
  pdftitle={The Heawood graph and the indicated domination game},
  pdfauthor={}}

\begin{document}
\begin{frontmatter}
\titledata{The Heawood graph and the indicated domination game}{}


\authordata{Xiaoxiang Yu}
  {School of Mathematics and Statistics, Beijing Jiaotong University, Beijing 100044, P. R. China}
  {yuxiaoxiang1126@163.com }{}

\authordata{Federico Romaniello}{Department for Humanistic, Scientific and Social Innovation, University of Basilicata, Potenza, 85100, Italy}{federico.romaniello@unibas.it}{} 

\authordatatwo
  {Jianbing Liu}{jbliu1@bjtu.edu.cn}{}
  {Rong-Xia Hao}{rxhao@bjtu.edu.cn}{}
  {School of Mathematics and Statistics, Beijing Jiaotong University,
   Beijing 100044, P. R. China}
\makeatletter
\advance\@authorcount by 2\relax
\makeatother

\begin{abstract}
We study a problem of Brešar et al. asking whether the indicated domination number equals the independence number for every bipartite graph. We show that the Heawood graph is a counterexample of minimum order among connected cubic bipartite graphs, having independence number seven and indicated and Grundy domination numbers both equal to eight, and we characterise it by a two-turn opening property. Using assemblies of vertex-deleted Heawood graphs, we construct connected cubic bipartite graphs where the difference between these parameters is arbitrarily large. Finally, we establish an improved linear upper bound for cubic graphs and bound the indicated domination number for the incidence graphs of finite projective planes.
\end{abstract}
\end{frontmatter}

\section{Introduction}\label{sec:introduction}

All graphs considered here are finite, non-empty and, unless stated
otherwise, simple. For graph-theoretic terminology and notation not
defined here, we refer to Bondy and Murty~\cite{BM2008}.

Let $G$ be a graph of order $n(G)=|V(G)|$. The open neighbourhood
$N_G(v)$ of a vertex $v$ consists of its neighbours, and its closed
neighbourhood is $N_G[v]=N_G(v)\cup\{v\}$. For $S\subseteq V(G)$, we
write $N_G[S]=\bigcup_{v\in S}N_G[v]$, omitting the subscript when the
graph is clear. The vertices dominated by $S$ are those in $N[S]$.
The set $S$ is \emph{dominating} if $N[S]=V(G)$.
The independence number $\alpha(G)$ is the maximum size of a set of
pairwise nonadjacent vertices. The vertex cover number $\tau(G)$ is the
minimum size of a set meeting every edge. Taking complements gives
$\tau(G)=n(G)-\alpha(G)$.

Bre\v{s}ar, Bujt\'as, Ir\v{s}i\v{c}, Rall and Tuza~\cite{BBIRT2024}
introduced the \emph{indicated domination game} as a variant of the
domination game of Bre\v{s}ar, Klav\v{z}ar and Rall~\cite{BKR2010}.
Initially no vertex is selected. If $S$ is the set of vertices selected
so far, Dominator indicates a vertex $v$ outside $N[S]$, and Staller
selects a vertex of $N[v]$. The \emph{gain} of a turn is the number of
vertices it newly dominates. Play ends when $S$ is a dominating set.
Dominator minimises the number of turns and Staller maximises it.
Under optimal play, this number is the \emph{indicated domination number}
$\gi(G)$. Henning and Rall~\cite{HR2025} studied the corresponding game
for total domination.

The \emph{upper domination number} $\Gamma(G)$ is the maximum size of an inclusion-minimal
dominating set. Given such a set $M$, Staller can choose a vertex of $M$
in response to every indication. By minimality, play then continues
until all of $M$ has been selected. In particular, Staller can play
throughout in a maximum independent set, since every such set is an
inclusion-minimal dominating set.

A \emph{legal dominating sequence} is a sequence of distinct vertices in which each
term dominates at least one new vertex and the whole sequence dominates
$G$. Its maximum length is the \emph{Grundy domination number}
$\ggr(G)$~\cite{BGMRR2014}. Staller's selections form a legal dominating
sequence, so $\Gamma(G)\leq\gi(G)\leq\ggr(G)$~\cite{BBIRT2024}.

Cockayne et al.~\cite{CFPT1981} proved that $\Gamma(G)=\alpha(G)$ for
bipartite graphs. Bre\v{s}ar et al.~\cite{BBIRT2024} proved that the
indicated domination number equals the independence number for trees
and grids. In Question~7.2 they asked whether this equality holds for
every bipartite graph. In a cubic bipartite graph of order $n$, the
two parts have equal size, and Hall's theorem gives a perfect matching.
Each part is therefore a maximum independent set of size $n/2$. The
same argument applies to cubic bipartite multigraphs. Bre\v{s}ar et
al.\ proved that the indicated domination number of every connected
cubic bipartite graph of order at most twelve is half its order. In
Question~7.3 of~\cite{BBIRT2024}, they asked whether this holds for all
cubic bipartite graphs.

The Heawood graph $\HH$, the incidence graph of the Fano plane, answers
both questions. In Theorem~\ref{prop:heawood} we prove that
$\gi(\HH)=8$ and $\alpha(\HH)=7$, so $\HH$ has minimum order among
connected cubic bipartite counterexamples. In the first two turns
Staller selects a nonincident point and line. From this position, every
indication admits a reply that dominates just one new vertex.
More generally, suppose that Staller can force such a position in two
turns on a connected $d$-regular bipartite graph, where $d\geq2$.
Theorem~\ref{thm:opening} gives an upper bound of $4d+2$ on the order
of the graph, with equality precisely for $\HH$.

The indicated domination number is additive over components, since
Dominator can finish one component at a time and Staller can play
optimally in each component regardless of the order of indications.
Disjoint copies of $\HH$ thus give graphs with arbitrarily large
difference $\gi-\alpha$. To obtain connected examples, we replace each
vertex of a loopless cubic multigraph $B$ with a copy of $\HH$ from which one
vertex has been deleted. The three edges at each vertex of $B$ are
attached to distinct degree-two vertices of its copy. If $B$ has order
$b$, Theorem~\ref{thm:assembly} gives $7b$ as the indicated domination
number of the resulting graph. By Proposition~\ref{prop:independence},
this exceeds its independence number by the vertex cover number of $B$.
If $B$ is connected and bipartite, so is the resulting graph, and the
difference is $b/2$.

Bre\v{s}ar et al.~\cite{BBIRT2024} also asked whether there is a constant
$C<1$ such that $\gi(G)\leq Cn(G)$ for every connected cubic graph $G$.
In Section~\ref{sec:bounds} we prove that $C=20/33$ suffices for all
cubic graphs and give a sharper bound for connected graphs of order
greater than $22$.

The interpretation of $\HH$ as the incidence graph of the Fano plane
suggests asking how the game depends on the order of the plane. In
Section~\ref{sec:planes} we show that if a finite projective plane has
$v$ points, then its incidence graph has Grundy domination number
$v+1$ and indicated domination number in $\{v,v+1\}$.
Section~\ref{sec:remarks} collects open problems on cubic bipartite
graphs and incidence graphs of projective planes.

\section{The Heawood graph and its assemblies}\label{sec:main}

We identify the vertices of $\HH$ with the seven points and seven lines of
the Fano plane, adjacency being incidence. Each line contains three points,
and each point lies on three lines. Write $pq$ for the unique line
through distinct points $p$ and $q$. Any two distinct lines meet in
a unique point. The Fano plane is self-dual, so $\HH$ has an automorphism
interchanging its point and line vertices.

\begin{lemma}\label{lem:four-vertices}
Every set of four vertices of $\HH$ dominates at least ten vertices.
\end{lemma}
\begin{proof}
By duality, we may assume that at least two of the four vertices are
points. If two are points and two are lines, the two points meet five
lines and the two lines contain five points. If three are points, a
collinear triple meets all seven lines. A noncollinear triple meets six
lines, and the selected line contains a fourth point. Finally, if all
four vertices are points, at most one line avoids them, since such a
line must consist of the three remaining points. In each case at least
ten vertices are dominated.
\end{proof}

\begin{figure}[!b]
\centering
\begin{tikzpicture}[every node/.style={font=\small},line cap=round,
  point mark/.style={circle,draw=black,inner sep=0pt,minimum size=5.8pt},
  line mark/.style={rectangle,draw=black,inner sep=0pt,minimum size=5.4pt},
  active/.style={fill=white,line width=.8pt},
  dominated/.style={fill=black!45,line width=.5pt},
  chosen/.style={fill=black,line width=.7pt},
  node label/.style={inner sep=1pt}]
  \begin{scope}[x=1.14cm,y=1.14cm,yshift=-.47cm]
    \coordinate (p) at (0,0);
    \coordinate (a) at (90:2);
    \coordinate (b) at (210:2);
    \coordinate (c) at (330:2);
    \coordinate (ab) at (150:1);
    \coordinate (bc) at (270:1);
    \coordinate (ca) at (30:1);
    \draw[dash pattern=on 2.7pt off 1.8pt,line width=.7pt] (a)--(b)--(c)--cycle;
    \draw[black!65,line width=.65pt] (a)--(bc) (b)--(ca) (c)--(ab);
    \draw[line width=1.05pt,preaction={draw=white,line width=2.5pt}] (p) circle[radius=1];
    \foreach \v in {a,b,c} \node[point mark,active] at (\v) {};
    \foreach \v in {ab,bc,ca} \node[point mark,dominated] at (\v) {};
    \node[point mark,chosen] at (p) {};
    \node[node label,above=3pt] at (a) {$q$};
    \node[node label] at (.27,.04) {$p$};
    \node[node label,anchor=east] at (-.08,1.28) {$pq$};
    \node[node label,above right=2pt] at (ca) {$x$};
    \node[node label,right=1pt] at (1.4,-.36) {$m$};
    \node[node label] at (120:.76) {$\ell$};
  \end{scope}
  \node at (0,-2.47) {Fano plane};

  \begin{scope}[xshift=5.65cm]
    \foreach \id [count=\i from 0] in
      {P,LB,CApt,CAline,C,BCline,B,ABline,A,LA,BCpt,E,ABpt,LC}
      \coordinate (g-\id) at ({90+360*(\i-5.5)/14}:1.74);
    \foreach \lin/\pts in {E/{ABpt,BCpt,CApt},LA/{P,A,BCpt},
      LB/{P,B,CApt},LC/{P,C,ABpt},ABline/{A,ABpt,B},
      BCline/{B,BCpt,C},CAline/{C,CApt,A}}
      \foreach \pt in \pts
        \draw[black!55,line width=.52pt] (g-\lin)--(g-\pt);
    \draw[line width=1.55pt] (g-A)--(g-LA) (g-CAline)--(g-CApt);
    \foreach \id in {A,B,C} \node[point mark,active] at (g-\id) {};
    \foreach \id in {ABline,BCline,CAline} \node[line mark,active] at (g-\id) {};
    \foreach \id in {ABpt,BCpt,CApt} \node[point mark,dominated] at (g-\id) {};
    \foreach \id in {LA,LB,LC} \node[line mark,dominated] at (g-\id) {};
    \node[point mark,chosen] at (g-P) {};
    \node[line mark,chosen] at (g-E) {};
    \node[node label,above left=2pt] at (g-A) {$q$};
    \node[node label,left=3pt] at (g-LA) {$pq$};
    \node[node label,above right=2pt] at (g-CAline) {$m$};
    \node[node label,right=3pt] at (g-CApt) {$x$};
    \node[node label,below right=2pt] at (g-P) {$p$};
    \node[node label,below left=2pt] at (g-E) {$\ell$};
    \node at (0,-2.47) {Heawood graph};
  \end{scope}
\end{tikzpicture}
\caption{The position after Staller selects $p$ and $\ell$, with $\ell$
drawn as a circle in the Fano plane. The three open corners and three
dashed sides correspond to the six undominated vertices of $\HH$.
In the graph, circles denote points and squares denote lines. Selected
vertices are black, their neighbours grey and undominated vertices white.
The bold edges join $q$ to $pq$ and $m$ to $x=m\cap\ell$, pairing each
of these indications with Staller's reply.}
\label{fig:heawood-strategy}
\end{figure}
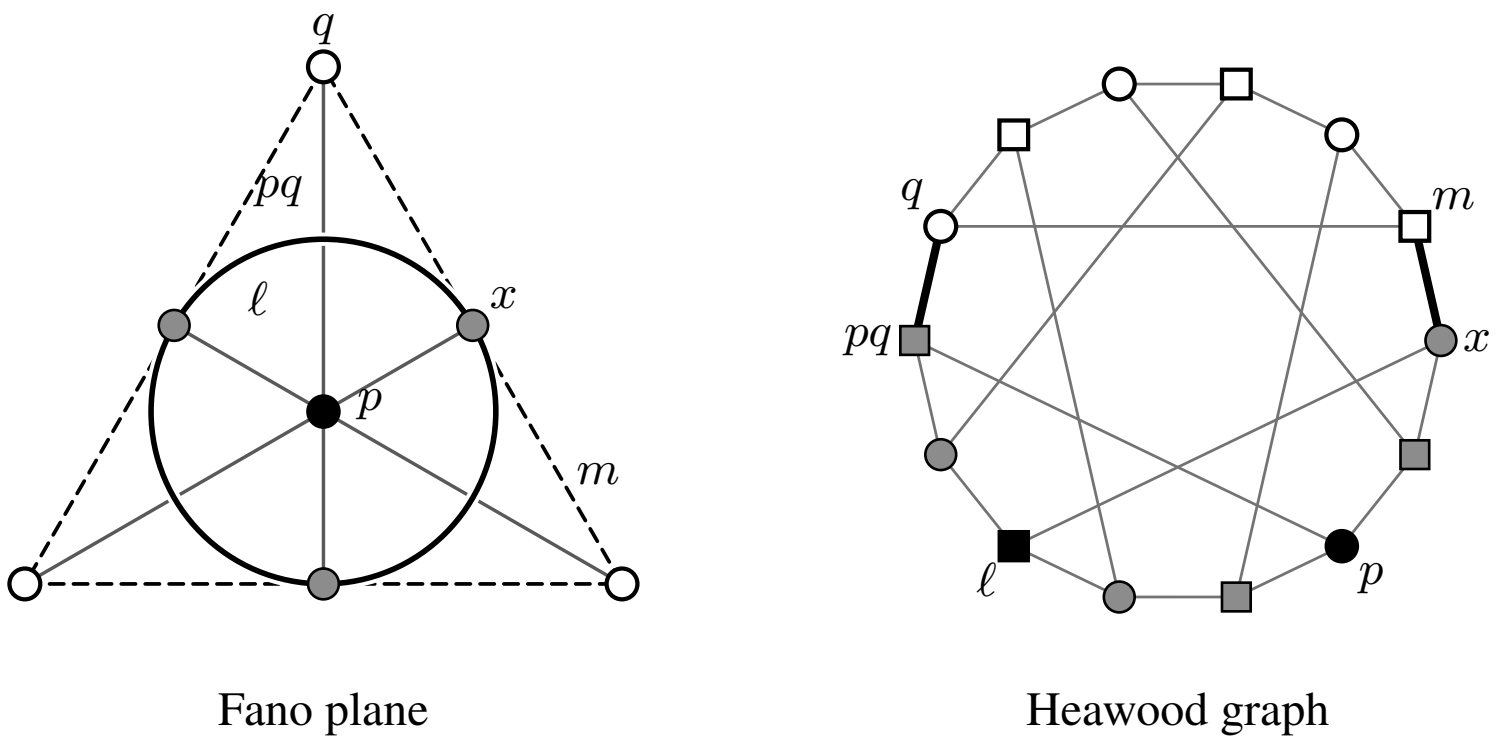

\begin{theorem}\label{prop:heawood}
The Heawood graph $\HH$ satisfies $\gi(\HH)=\ggr(\HH)=8>\alpha(\HH)=7$.
It has minimum order among connected cubic bipartite graphs $G$ with
$\gi(G)>\alpha(G)$.
\end{theorem}
\begin{proof}
By point--line duality, we may assume that Dominator first indicates a
point $p$. Staller selects $p$ itself and, on the second turn, a line
$\ell$ avoiding $p$. Indeed, if Dominator indicates a line, that line
already avoids $p$. If Dominator indicates a point $q$, either line
through $q$ other than $pq$ is a legal choice.

Since $p$ and $\ell$ are nonincident, $N[p]$ and $N[\ell]$ are disjoint
sets of four vertices. Let $U$ be the set of six vertices left
undominated, as in Figure~\ref{fig:heawood-strategy}. For a point $q\in U$,
the line $pq$ is already dominated, and its three points are $p$, $q$
and $pq\cap\ell$. Hence $N[pq]\cap U=\{q\}$. Dually, for a line
$m\in U$, we have $N[m\cap\ell]\cap U=\{m\}$. Whenever Dominator
indicates a vertex of $U$, the corresponding reply therefore dominates
only that vertex for the first time. This forces six further turns, so
$\gi(\HH)\geq8$.

In any legal dominating sequence of length at least four, the first four
selections dominate at least ten vertices by Lemma~\ref{lem:four-vertices}.
Each subsequent selection must dominate a new vertex, so there are at
most four more. Thus $8\leq\gi(\HH)\leq\ggr(\HH)\leq8$.

Since $\HH$ is cubic and bipartite of order fourteen, $\alpha(\HH)=7$.
Any smaller connected cubic bipartite graph $G$ has even order at most
twelve and hence satisfies $\gi(G)=\alpha(G)$ by Proposition~5.5
of~\cite{BBIRT2024}.
\end{proof}

After the first two turns in Figure~\ref{fig:heawood-strategy}, every
undominated vertex admits a reply of gain one. Such a reply remains
available as long as that vertex is undominated. The following lemma
bounds the number of undominated vertices in a position of this kind.

\begin{lemma}\label{lem:two-gain}
Let $S$ be the set of selected vertices in a position on a $d$-regular
graph $G$, where $d\geq2$, and put $D=N[S]$ and $U=V(G)\setminus D$.
If every vertex of $U$ admits a reply of gain one, then
$|S|+2|U|\leq n(G)$. At equality, the edges between $U$ and $D\setminus S$
form a perfect matching between these two sets.
\end{lemma}
\begin{proof}
Let $I$ be the isolated vertices of $G[U]$. Each $v\in U\setminus I$
has a response $x_v$ with $N[x_v]\cap U=\{v\}$. It lies outside
$U$. Selecting $v$ would also dominate an undominated neighbour,
and selecting any other vertex of $U$ would dominate that vertex as
well as $v$. Hence $x_v\in D\setminus S$. These witnesses are
distinct and none is adjacent to $I$.

All neighbours of $I$ lie in $D\setminus S$. Each of them has a
neighbour in $S$ and hence at most $d-1$ neighbours in $I$. If $I$
is nonempty, its $d|I|$ incident edges therefore reach more than
$|I|$ neighbours. Counting these together with the witnesses gives
$|D\setminus S|\geq|U|$, or equivalently $|S|+2|U|\leq n(G)$.
At equality, $I$ is empty and the witnesses exhaust $D\setminus S$.
Each $x_v$ has $v$ as its only neighbour in $U$, giving the matching.
\end{proof}

We now apply the lemma to a position reached after two turns.

\begin{theorem}\label{thm:opening}
Let $G$ be a connected $d$-regular bipartite graph of order $n$,
where $d\geq2$. Suppose that Staller can respond to the first two indications,
however Dominator plays, so that every remaining undominated vertex
admits a response of gain one.
Then $n\leq 4d+2$, with equality if and only if $G$ is the Heawood graph.
\end{theorem}
\begin{proof}

Fix a possible first selection $p$ of Staller's strategy. With $p$
fixed, every second indication admits a reply $y$ that leaves each
undominated vertex with a response of gain one. Put $S=\{p,y\}$,
$D=N[S]$ and $U=V(G)\setminus D$.
Since $|N[p]\cup N[y]|\leq2d+2$, Lemma~\ref{lem:two-gain} gives
\[
 2+2(n-2d-2)\leq 2+2|U|\leq n.
\]
Suppose that $n=4d+2$. The two closed neighbourhoods are disjoint,
and $|U|=2d$. Equality in Lemma~\ref{lem:two-gain} gives each vertex
of $U$ exactly one neighbour in $D$, and each vertex of
$D\setminus S$ exactly one in $U$. Hence $G[U]$ is $(d-1)$-regular.
Its two parts have equal size, as do those
of $G[D]$. Since $N[p]$ and $N[y]$ are disjoint, this forces $p$ and
$y$ to belong to opposite parts.

Call the vertices in the part of $p$ points and those in the other
part lines. Each part has $2d+1$ vertices. For any
line $L$ avoiding $p$, Dominator may indicate $L$ on the second turn.
The only response on the opposite side from $p$ is $L$ itself.
Thus the preceding conclusions hold for \emph{every} such $L$.
A line $M\ne L$ avoiding $p$ is in $U$ and meets $L$ in precisely one
point, its sole neighbour in $D$. A line $M$ through $p$ is in
$D\setminus S$ and meets $L$ in $d-2$ points, since it contains $p$
and exactly one point in $U$.

Fix a line $M$ through $p$. Counting at its $d$ points, its
intersections with all other lines have total multiplicity $d(d-1)$.
The $d+1$ lines avoiding $p$ account for $(d+1)(d-2)$, leaving
exactly two incidences for the other $d-1$ lines through $p$.
Each contains $p$, so $d\leq3$. If $d=2$, the two lines through $p$
have both points in common and give a $4$-cycle component, contrary
to connectedness at order ten. If $d=3$, the other two lines through
$p$ each meet $M$ only at $p$. Thus all pairs of lines meet exactly
once. The seven triples contain $7\binom32=\binom72$ distinct pairs
of points, giving the Fano plane and hence its incidence graph $\HH$.
The converse is the strategy in Theorem~\ref{prop:heawood}.
\end{proof}

By additivity, $k$ disjoint copies of $\HH$ have indicated
domination number $8k$ and independence number $7k$. To obtain a growing
difference in connected graphs, we delete one vertex from each copy.

\subsection{Playing on a Heawood fragment}\label{sec:fragment}

Fix a point $p$ of the Fano plane and let $T=\HH-p$. The three
line vertices through $p$ have degree two in $T$. We call them the
\emph{ports} of $T$ and the remaining ten vertices its \emph{core}.
A \emph{block} is a copy of $T$.

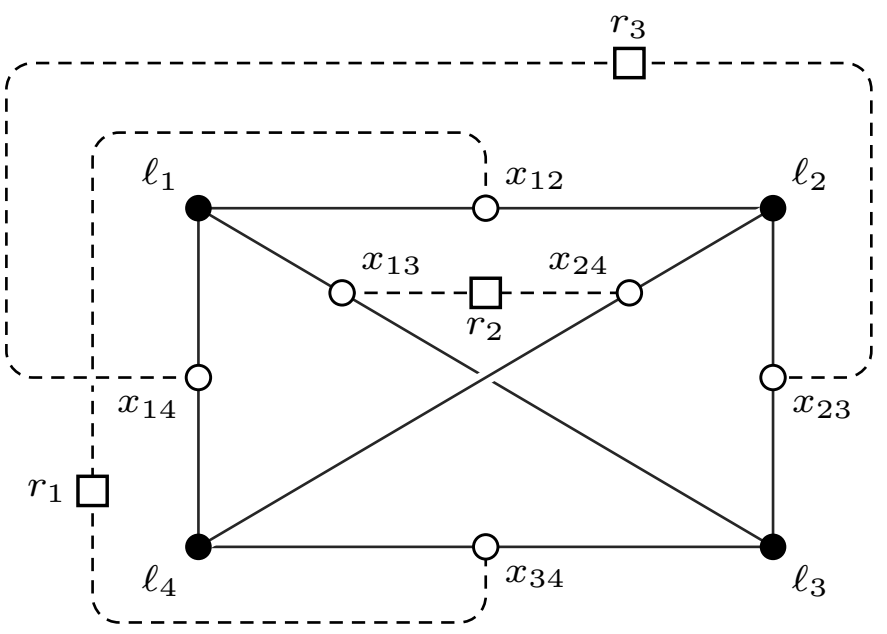
\begin{figure}[htb]
\centering
\begin{tikzpicture}[
  every node/.style={font=\scriptsize},
  every label/.style={inner sep=1.5pt},
  original/.style={circle,draw=black,fill=black,inner sep=0pt,minimum size=4.6pt},
  point/.style={circle,draw=black,fill=white,inner sep=0pt,minimum size=4.6pt,line width=.55pt},
  port/.style={rectangle,draw=black,fill=white,inner sep=0pt,minimum size=5.5pt,line width=.65pt},
  core edge/.style={draw=black!85,line width=.5pt},
  port edge/.style={draw=black,line width=.5pt,dash pattern=on 2.5pt off 1.7pt,rounded corners=5pt}]
  \coordinate (l1) at (-1.9,1.12);
  \coordinate (l2) at ( 1.9,1.12);
  \coordinate (l3) at ( 1.9,-1.12);
  \coordinate (l4) at (-1.9,-1.12);
  \coordinate (x12) at (0,1.12);
  \coordinate (x34) at (0,-1.12);
  \coordinate (x14) at (-1.9,0);
  \coordinate (x23) at ( 1.9,0);
  \coordinate (x13) at (-.95,.56);
  \coordinate (x24) at ( .95,.56);
  \coordinate (r1) at (-2.6,-.75);
  \coordinate (r2) at (0,.56);
  \coordinate (r3) at (.95,2.08);
  \draw[core edge] (l1)--(x12)--(l2)--(x23)--(l3)--(x34)--(l4)--(x14)--cycle;
  \draw[core edge] (l1)--(x13)--(l3);
  \draw[core edge,preaction={draw=white,line width=2.8pt}] (l4)--(x24)--(l2);

  \draw[port edge] (x12)--(0,1.62)--(-2.6,1.62)--(r1);
  \draw[port edge] (r1)--(-2.6,-1.62)--(0,-1.62)--(x34);
  \draw[port edge] (x13)--(r2)--(x24);
  \draw[white,line width=3pt] (-2.69,0)--(-2.51,0);
  \draw[port edge] (x14)--(-3.17,0)--(-3.17,2.08)--(r3);
  \draw[port edge] (r3)--(2.55,2.08)--(2.55,0)--(x23);

  \foreach \i/\pos in {1/above left,2/above right,3/below right,4/below left}
    \node[original,label=\pos:{$\ell_{\i}$}] at (l\i) {};
  \foreach \ij/\pos in {12/above right,13/above right,14/below left,
                         23/below right,24/above left,34/below right}
    \node[point,label=\pos:{$x_{\ij}$}] at (x\ij) {};
  \foreach \i/\pos in {1/left,2/below,3/above}
    \node[port,label=\pos:{$r_{\i}$}] at (r\i) {};
\end{tikzpicture}
\caption{The fragment $T=\HH-p$ as a subdivided $K_4$ with three ports.
The filled circles are the four lines $\ell_i$ avoiding $p$, and the open
circles are their intersections $x_{ij}=\ell_i\cap\ell_j$. The squares
represent the three lines through $p$. Dashed edges join each port to
the subdivision vertices on a pair of opposite edges of $K_4$.}
\label{fig:fragment}
\end{figure}

Figure~\ref{fig:fragment} shows how to read the incidences of $T$ on a
subdivided $K_4$. Its four original vertices are the lines avoiding $p$.
Each point other than $p$ is the intersection of two such lines and
subdivides the edge joining them. The two points on a line through $p$
belong to disjoint pairs of these four lines. Each port therefore joins
the subdivision vertices on a pair of opposite edges.

An edge joining two blocks has a port at either end. Staller can therefore
play inside a block even when some of its ports are dominated from outside.

\begin{lemma}\label{lem:lower}
Suppose that all core vertices of $T$ are initially undominated.
Ports may already be dominated, and play outside $T$ may dominate
further ports, provided it leaves the core unaffected. By responding
inside $T$, Staller can force at
least seven selections in $T$ before its core is dominated, even if
turns in $T$ are interspersed with turns elsewhere.
\end{lemma}
\begin{proof}
Until a core vertex is indicated, Staller responds to every indicated
port by selecting that port itself. Let $a$ be the number of ports
selected before the first core indication, where $0\leq a\leq3$.
These ports dominate $2a$ distinct points. The four lines avoiding
$p$ are still undominated, so a core indication must eventually occur.

When the first core vertex is indicated, Staller selects a line
$\ell$ avoiding $p$. For an indicated line, the reply is the line
itself. For an indicated point, either nonport line through it will do.
Each of the $a$ previously selected ports contains one point on $\ell$
and one outside $\ell$. The undominated core therefore consists of
$3-a$ points outside $\ell$ and three nonport lines other than $\ell$.

For every subsequent indication of a core point $q$, Staller selects
the port $pq$. Its other point lies on $\ell$ and is already dominated.
For an indicated nonport line $m$, the reply is $m\cap\ell$. This point
is already dominated, and its incident lines are $m$, $\ell$ and a
port. Hence $m$ is the only new core vertex. If a port is
indicated, Staller selects that port itself. Since every port has
exactly one point outside $\ell$, this response also dominates at
most one new core vertex.

Thus each further selection inside the block dominates at most one
of the $6-a$ remaining core vertices. External domination of ports
does not change this conclusion. At least
$a+1+(6-a)=7$ selections inside $T$ are necessary.
\end{proof}

When a port is indicated, Staller may select its external neighbour.
In the local game on $T$, we count this as one turn which dominates
only the indicated port.

\begin{lemma}\label{lem:upper}
Suppose that the four original vertices in the $K_4$ description of $T$
are initially undominated. Any subset of the subdivision vertices and
ports may already be dominated. Dominator can dominate $T$ in at most
seven turns, including turns on which Staller selects the external
neighbour of an indicated port.
\end{lemma}
\begin{proof}
Dominator first indicates undominated original vertices until all four
are dominated. This takes at most four turns. On each turn Staller
selects either the indicated vertex or a subdivision vertex on an
incident edge. Once this phase ends, Dominator can finish in at most
one turn for each vertex still undominated, since even an external
response dominates the indicated port.

Let $r$ be the number of original vertices selected by Staller during
the first phase.

\textbf{Case 1.} $r\geq2$. If at least three original vertices were
selected, every edge of $K_4$ has a selected endpoint. Hence all
subdivision vertices are dominated and only the three ports may
remain. If exactly two original vertices were selected, all subdivision
vertices are dominated except possibly the one on the edge joining
the two unselected originals. To dominate those two
originals, Staller must have selected at least one subdivision vertex,
thereby dominating its port. Thus in either case at most three vertices
remain, and the total is at most $4+3=7$.

\textbf{Case 2.} $r=1$. Let $u$ be the selected original vertex.
Only the three ports and the three subdivision vertices on edges
avoiding $u$ can remain undominated. Each subdivision
vertex selected during the first phase accounts for a distinct member
of this set. If its edge avoids $u$, count the subdivision vertex
itself. If its edge meets $u$, count its port. These ports are
distinct because the three edges at $u$ belong to different pairs of
opposite edges. Thus, if Staller selected $s$ subdivision vertices,
the first phase took $1+s$ turns and at most $6-s$ vertices remain.
Again the total is at most seven.

\textbf{Case 3.} $r=0$. All $s$ selections in the first phase were
subdivision vertices. If at least two ports are dominated, at most
$6-s$ subdivision vertices and one port remain, giving at most seven
turns in all.

It remains to consider the case in which only one port is dominated.
Every selected subdivision vertex then belongs to the pair of
opposite edges corresponding to that port. Both subdivision vertices
must have been selected to dominate all four originals. After these
two turns, the only vertices that may remain undominated are the four
subdivision vertices on the complementary $4$-cycle and the other
two ports.

If at most five of these six vertices remain undominated, the bound
follows. Otherwise, Dominator indicates any subdivision vertex $x$ on
the cycle. Selecting $x$ or its port dominates both of them. Selecting
either original endpoint of the edge containing $x$ dominates the two
subdivision vertices incident with that endpoint on the cycle. Every
response therefore dominates at least two new vertices. After the
third turn at most four vertices remain, so seven turns suffice.
\end{proof}

\subsection{Assembling the fragments}\label{sec:composition}

Staller can now play in the block containing each indication, while
Dominator can finish the blocks one at a time.

\begin{theorem}\label{thm:assembly}
Let $t\geq2$ be even. Take $t$ disjoint copies of $T$ and add a perfect
matching on their $3t$ ports, with each added edge joining distinct
blocks. The resulting graph $G$ is simple and cubic, and
\begin{linenomath*}
\[
 \gi(G)=7t.
\]
\end{linenomath*}
It is connected if the graph recording the connections between blocks is connected.
\end{theorem}
\begin{proof}
Each port acquires one neighbour in a different block, so $G$ is simple
and cubic. Since each block is connected, the stated condition gives
connectedness.

Staller answers each indication within its block, using
Lemma~\ref{lem:lower}. A selection outside that block can dominate only a
port there. Its core therefore requires at least seven local selections,
giving $\gi(G)\geq7t$.

For the reverse bound, Dominator finishes one block at a time. Before
the first indication in a block, any vertex already selected there was
a port answering an external indication. Its local neighbours are all
subdivision vertices, so all four original vertices are still undominated.
While Dominator plays in the block, an external response can dominate
only the indicated port there. Lemma~\ref{lem:upper} now gives at most
seven turns per block, counted where each indication is made.
\end{proof}

For a loopless cubic multigraph $B$, let $\TT(B)$ be obtained by
replacing every vertex with a copy of $T$, assigning its three incident
edges to distinct ports, and joining the assigned ports along the
edges of $B$. This notation includes any choice of the assignments.
Parallel edges are allowed in the base, but $\TT(B)$ is simple.

\begin{proposition}\label{prop:independence}
For every loopless cubic multigraph $B$,
\begin{linenomath*}
\[
 \alpha(\TT(B))=6n(B)+\alpha(B),\qquad
 \gi(\TT(B))-\alpha(\TT(B))=\tau(B).
\]
\end{linenomath*}
\end{proposition}
\begin{proof}
An independent set of size seven in $\HH$ meets every edge, so its
complement is independent as well. The bipartition of the connected
graph $\HH$ is unique up to interchange. Thus the only independent set
of size seven in $T$ consists of all seven lines, including the ports.

For an independent set of $\TT(B)$, the blocks contributing seven vertices
therefore form an independent set of $B$. Every other block contributes at
most six, giving the upper bound $6n(B)+\alpha(B)$. To attain it, take
all seven lines in the blocks of a maximum independent set of $B$ and all
six points in the other blocks. The second equality follows from
Theorem~\ref{thm:assembly} and $\tau(B)=n(B)-\alpha(B)$.
\end{proof}

Suppose now that $B$ is cubic and bipartite. Choose the bipartitions of
the blocks so that ports of adjacent blocks belong to opposite
classes. Then $\TT(B)$ is cubic and bipartite, and is connected
whenever $B$ is connected. Theorem~\ref{thm:assembly} and
Proposition~\ref{prop:independence} therefore give
\begin{linenomath*}
\begin{equation}\label{eq:bipartite-values}
 n(\TT(B))=13n(B),\qquad \gi(\TT(B))=7n(B),\qquad
 \alpha(\TT(B))=\frac{13}{2}n(B).
\end{equation}
\end{linenomath*}

For $k\geq1$, let $M_k$ have two edges from $A_j$ to $B_j$ and one from
$A_j$ to $B_{j-1}$ for each $j$ modulo $k$. This connected cubic bipartite
base has order $2k$. By~\eqref{eq:bipartite-values}, $G_k=\TT(M_k)$ has
order $26k$, with $\gi(G_k)=14k$ and $\alpha(G_k)=13k$. Thus the difference $\gi(G_k)-\alpha(G_k)=k$ is unbounded among
connected cubic bipartite graphs. For $k=1$, the two copies of $T$ are
joined by a matching of their ports.

\section{Upper bounds for cubic graphs}\label{sec:bounds}

Bre\v{s}ar et al.~\cite{BBIRT2024} asked whether there is a constant
$C<1$ such that $\gi(G)\leq Cn(G)$ for every connected cubic graph.
Lemma~\ref{lem:two-gain} already gives $C=2/3$. Indeed, Dominator
indicates a vertex forcing a gain of at least two whenever possible.
If this strategy takes $t$ turns and leaves $u$ vertices undominated,
then $2t+u\leq n$ by the gains, and $t+2u\leq n$ by the lemma.
Dominator can finish in at most $u$ further
turns, so $\gi(G)\leq t+u\leq 2n/3$. This applies in every regular
graph of degree at least two.

In a cubic graph Dominator can do better by first forcing larger gains.

\begin{theorem}\label{prop:linear}
Every cubic graph $G$ of order $n$ satisfies $\gi(G)\leq 20n/33$.
If $G$ is connected, then $\gi(G)\leq (109n+2)/180$ as well.
\end{theorem}
\begin{proof}
Let $S$ be the set of vertices selected so far, initially empty.
Whenever a vertex lies outside all balls of radius three centred at
vertices of $S$, Dominator indicates such a vertex. Every possible
response has a closed neighbourhood disjoint from $N[S]$, so the
turn has gain four. When this stage ends, the balls of radius three
around the $s$ selected vertices cover $G$. Each ball contains at
most $1+3+6+12=22$ vertices, so $n\leq22s$.

If $G$ is connected, each indication after the first can be chosen at
distance four from $S$ until the balls cover $G$. The new centre is then
within distance five of an earlier centre, so their radius-three balls
share at least two vertices on a shortest path between them. Hence
$n\leq22+20(s-1)=20s+2$.

Dominator next forces a gain of at least three for as long as
possible. Let $a$ and $b$ count all turns of gain four and three
since the start. Then $a\geq s$ and $|D|=4a+3b$, where $D=N[S]$.
A response adjacent to a previous selection has gain at most two,
since it and that selection are already dominated. Thus $S$ is
independent and $G[D]$ contains at least $3(a+b)$ edges.

Write $e(D,U)$ for the number of edges from $D$ to the undominated
set $U$. We claim that $|U|\leq e(D,U)$. Assign to each vertex in
$U$ with a neighbour in $D$ one edge to $D$. Any remaining vertex
$v$ has all three neighbours in $U$. Indicating $v$ now admits a
response of gain at most two. This must be a neighbour $x$ of $v$
whose other two neighbours lie in $D$. Such an $x$ is distinct for
each $v$ and has a second edge to $D$, which we assign to $v$.
Summing degrees in $D$ gives
\[
 n\leq |D|+e(D,U)=4|D|-2|E(G[D])|\leq10a+6b.
\]
Consequently $6(2a+b)\geq n+2a\geq n+2s$.

Finally Dominator forces a gain of at least two for as long as
possible. Write $t$ for the total number of turns played by then and
$u$ for the number of undominated vertices. The earlier gains give
$2t+u+2a+b\leq n$, while Lemma~\ref{lem:two-gain} gives
$t+2u\leq n$. Adding and using the preceding estimate yields
\[
 \gi(G)\leq t+u\leq\frac{2n-(2a+b)}3
             \leq\frac{11n-2s}{18}.
\]
The two bounds follow from $s\geq n/22$ and, for connected $G$,
$s\geq(n-2)/20$.
\end{proof}

For connected graphs, the general bound is smaller when $n<22$, the two
agree at $n=22$, and the connected bound is smaller when $n>22$.

\section{Incidence graphs of projective planes}\label{sec:planes}

By Theorem~\ref{prop:heawood}, Staller can force eight turns on the
Heawood graph, one more than the number of points in the Fano plane. Can
Staller also force an extra turn on the incidence graph of a larger plane?

Let $\Pi$ be a projective plane of order $q\geq2$ with $v=q^2+q+1$
points. Its incidence graph $H_\Pi$ is $(q+1)$-regular and bipartite,
with $v$ points and $v$ lines. Hall's theorem gives a perfect matching,
so $\alpha(H_\Pi)=v$. By selecting only points, Staller can therefore
force at least $v$ turns.

The strict inequality $\ggr(H_\Pi)>\alpha(H_\Pi)$ follows from Bacs\'o et al.~\cite{BBKR2023}.
For $\mathrm{PG}(2,q)$, Gant's spectral calculation for zero
forcing~\cite{Gant2024} and Lin's rank bound~\cite{Lin2019} give
$\ggr(H_\Pi)\leq v+1$. We give a direct proof for every projective plane.

\begin{proposition}\label{prop:planes}
For every projective plane $\Pi$ of order $q\geq2$,
\begin{linenomath*}
\[
 q^2+q+1\leq\gi(H_\Pi)\leq\ggr(H_\Pi)=q^2+q+2.
\]
\end{linenomath*}
\end{proposition}
\begin{proof}
We use Lin's rank method~\cite{Lin2019}.
Let $A$ be the point--line incidence matrix. Then
$A^{\mathsf T}A=qI_v+J_v$, where $J_v$ is the all-one matrix.
The real matrix
\begin{linenomath*}
\[
 M=\begin{pmatrix}I_v&A\\ A^{\mathsf T}&qI_v\end{pmatrix}
\]
\end{linenomath*}
has rank $v+1$, since the Schur complement of $I_v$ is $-J_v$.
Its nonzero entries correspond exactly to equal or adjacent vertices.
Given a legal dominating sequence $x_1,\ldots,x_s$, choose a vertex
$y_i$ first dominated by $x_i$. These witnesses are distinct, and
$(M_{y_i,x_j})_{i,j=1}^s$ is triangular with nonzero diagonal.
Thus $s\leq v+1$.

For equality, fix a line $\ell$ and select all $q^2$ points outside
$\ell$. Each point dominates itself for the first time, and together
they dominate every line except $\ell$. Select a line $m\ne\ell$,
then the point $r=m\cap\ell$, and finally the remaining $q$ points
on $\ell$. These selections successively dominate a new point,
the line $\ell$, and the remaining points. The resulting legal
dominating sequence has length $q^2+1+1+q=v+1$.
The inequalities follow from $\alpha(H_\Pi)\leq\gi(H_\Pi)\leq
\ggr(H_\Pi)$.
\end{proof}

For $q\geq3$, we have $2(q^2+q+1)>4(q+1)+2$. By
Theorem~\ref{thm:opening}, Staller cannot then guarantee after two turns
that every remaining undominated vertex admits a reply of gain one.

\section{Open problems}\label{sec:remarks}

Let $\mathcal C$ be the class of connected cubic bipartite graphs. Put
\[
 \rho=\sup_{G\in\mathcal C}\frac{\gi(G)}{n(G)},\qquad
 \rho_\infty=\lim_{N\to\infty}
 \sup_{\substack{G\in\mathcal C\\ n(G)\geq N}}\frac{\gi(G)}{n(G)}.
\]
The Heawood graph has ratio $4/7$, and the graphs $G_k$ of
Section~\ref{sec:composition} have ratio $7/13$ at arbitrarily large
orders. Together with Theorem~\ref{prop:linear}, these examples give
\[
 \frac47\leq\rho\leq\frac{20}{33},\qquad
 \frac7{13}\leq\rho_\infty\leq\frac{109}{180}.
\]
Without the connectedness assumption, disjoint copies of $\HH$ give
ratio $4/7$ at arbitrarily large orders.

\begin{question}\label{q:ratios}
What are the values of $\rho$ and $\rho_\infty$?
\end{question}

In our connected assemblies, the three edges leaving a block form an
edge cut with a cycle on each side. A cubic graph is cyclically
$4$-edge-connected if every edge cut separating two subgraphs that
contain cycles has size at least four.

\begin{question}\label{q:cyclic}
Is $\gi(G)-n(G)/2$ unbounded among cyclically
$4$-edge-connected cubic bipartite graphs?
\end{question}

For a projective plane $\Pi$ of order $q$, Proposition~\ref{prop:planes}
gives two possible values for $\gi(H_\Pi)$. The Fano plane attains the
larger value. Write $H_q$ for the incidence graph of $\mathrm{PG}(2,q)$.

\begin{question}\label{q:planes}
For which projective planes $\Pi$ does $\gi(H_\Pi)=q^2+q+2$?
In particular, does $\gi(H_q)=q^2+q+1$ hold for every prime
power $q\geq3$?
\end{question}

\section*{Acknowledgments}

This research is partially supported by National Natural Science
Foundation of China (Nos.12471321, 12331013), Fundamental Research
Funds for the Central Universities (No. KYAJBRC26004536) and Hebei
Provincial Key Laboratory of Mathematical Theory and Analysis for
Network and Data Science.

The second author wishes to express his deep gratitude to the School of Mathematics and Statistics at Beijing Jiaotong University for hosting him and for their warm hospitality during an exchange visit under the Foreign Experts National Project, during which this manuscript was drafted.

\section*{Declaration of competing interest}

The authors declare that they have no known competing financial interests or personal
relationships that could have appeared to influence the work reported in this paper.

\end{document}